\documentclass[12pt]{amsart}
\usepackage[dvips]{graphicx}
\usepackage{amssymb,amsmath,amsthm}
\def\diG{\overrightarrow{G}}

\def\bold{\bf}

\def\0b{{\bold 0}}
\usepackage{bm}
\bmdefine{\Bzero}{0}
\bmdefine{\Bone}{1}
\def\Bone{{\bf 1}}
\def\RR{{\mathbb R}}
\def\ZZ{{\mathbb Z}}

\def\QQ{{\mathbb Q}}
\newtheorem{Theorem}{Theorem}[section]
\newtheorem{Lemma}[Theorem]{Lemma}
\newtheorem{Corollary}[Theorem]{Corollary}
\newtheorem{Proposition}[Theorem]{Proposition}
\newtheorem{Remark}[Theorem]{Remark}

\newtheorem{Example}[Theorem]{Example}

\newtheorem{Conjecture}[Theorem]{Conjecture}
\newtheorem{Question}[Theorem]{Question}
\begin{document}
\title{Centrally symmetric configurations of integer matrices}
\author{Hidefumi Ohsugi and Takayuki Hibi}
\thanks{
{\bf 2010 Mathematics Subject Classification:}
Primary 13P10; Secondary 52B20. \\
\hspace{5.5mm}{\bf Keywords:}
configuration, unimodular matrix, 
centrally symmetric configuration, 
Gr\"obner bases, 
unimodular triangulation,
Gorenstein ring.
}
\date{}
\begin{abstract}
The concept of centrally symmetric configurations
of integer matrices is introduced.
We study the problem when the toric ring of
a centrally symmetric configuration is normal
as well as is Gorenstein.  In addition, 
Gr\"obner bases of toric ideals of
centrally symmetric configurations will be discussed. 
Special attentions will be given to centrally
symmetric configurations of unimodular matrices and those of
incidence matrices of finite graphs.  
\end{abstract}
\maketitle
\section*{Introduction}
A {\em configuration} of $\RR^d$ is a matrix $A \in \ZZ^{d \times n}$, 
where $n = 1, 2, \ldots,$ for which
there exists a hyperplane ${\mathcal H} \subset \RR^d$ 
not passing the origin of $\RR^d$ such that each column vector 
of $A$ lies on ${\mathcal H}$.
Let $K$ be a field and 
$K[{\bf t}, {\bf t}^{-1}] 
= K[t_1, t_1^{-1}, \ldots, t_d, t_d^{-1}]$
the Laurent polynomial ring in $d$ variables over $K$.
Each column vector 
${\bf a} = [a_1, \ldots, a_d]^\top \in \ZZ^d$ 
$( = \ZZ^{d \times 1})$, where
$[a_1, \ldots, a_d]^\top$ is the transpose of
$[a_1, \ldots, a_d]$,
yields the Laurent monomial ${\bf t}^{\bf a} =
t_1^{a_1} \cdots t_d^{a_d}$.
Let $A \in \ZZ^{d \times n}$ be a configuration of $\RR^d$
with ${\bf a}_1, \ldots, {\bf a}_n$ its column vectors.   
The {\em toric ring} of $A$ 
is the subalgebra $K[A]$ of $K[{\bf t}, {\bf t}^{-1}]$ 
which is generated by the Laurent monomials
${\bf t}^{{\bf a}_1}, \ldots, {\bf t}^{{\bf a}_n}$.
Let $K[{\bf x}] = K[x_1, \ldots, x_n]$ be the polynomial
ring in $n$ variables over $K$ and define the surjective
ring homomorphism $\pi : K[{\bf x}] \rightarrow K[A]$
by setting $\pi(x_i) = {\bf t}^{{\bf a}_i}$ 
for $i = 1, \ldots, n$. 
We say that the kernel $I_A \subset K[{\bf x}]$ of $\pi$ is 
the {\em toric ideal} of $A$.
Finally, we write ${\rm Conv}(A)$ for the convex hull
of ${\bf a}_1, \ldots, {\bf a}_n$ in $\RR^d$.
Thus ${\rm Conv}(A) \subset \RR^d$ is 
an integral convex polytope, i.e.,
a convex polytope all of whose vertices have integer coordinates.
 
Given a matrix $A \in \ZZ^{d \times n}$, 
which is not necessarily a configuration of $\RR^d$, 
we introduce 
the {\it centrally symmetric configuration} $A^\pm \in \ZZ^{(d+1) \times (2n + 1)}$
of $\RR^{d + 1}$ as follows:
\[
A^{\pm} = 
\left[
\begin{array}{c|ccc|ccc}
0      &   &   &    & &   & \\
\vdots &   & A &    & &-A & \\
0      &   &   &    & &   & \\
\hline
1      & 1 & \cdots & 1 & 1 & \cdots & 1
\end{array}
\right].
\]

In the present paper we establish fundamental results 
on centrally symmetric configurations.  Especially we study centrally 
symmetric configurations of unimodular matrices and those of
incidence matrices of finite graphs.  
First, in Section $1$, by means of the notion of 
Hermite normal form of an integer matrix, we pay attention to 
the fact that  
the index $[\ZZ^d : \ZZ A]$, where $\ZZ A$ is
the abelian subgroup of $\ZZ^d$ generated by the column vectors
of an integer matrix $A \in \ZZ^{d \times n}$ of rank $d$,
is equal to the greatest common divisor of maximal minors of $A$.
Moreover, if $A \in \ZZ^{d \times n}$ is of rank $d$, then
there exists a nonsingular matrix $B$ such that $A' = B^{-1}A$
is an integer matrix satisfying $\ZZ A' = \ZZ^d$.
Second, in Section $2$, we study the centrally symmetric 
configuration of a unimodular matrix.  Two fundamental theorems
will be given.  Theorem \ref{unimodularinitial} says that
if $A \in \ZZ^{d \times n}$ is a unimodular matrix, then
there exists a reverse lexicographic order $<$ 
on $K[{\bf x}]$ such that
${\rm in}_< (I_{A^\pm})$ is a radical ideal.
Thus in particular $K[A^\pm]$ is normal.
Moreover, Theorem
\ref{unimodularGorenstein} guarantees that
the toric ring $K[A^\pm]$ of the centrally symmetric configuration 
$A^{\pm}$ of a unimodular matrix is Gorenstein.

On the other hand, we devote Sections $3$ and $4$ to the study
on toric rings and toric ideals of centrally symmetric configurations
of incidence matrices of finite graphs.     
Let $G$ be a finite connected graph on the vertex set
$[d] = \{ 1, \ldots, d \}$ and suppose that $G$ possesses
no loop and no multiple edge.
Let $E(G) = \{ e_1, \ldots, e_n \}$ denote
the edge set of $G$.  Let ${\bf e}_1, \ldots, {\bf e}_d$
stand for the canonical unit coordinate vector of $\RR^d$.
If $e = \{ i, j \}$ is an edge of $G$ with $i < j$, then
the column vectors $\rho(e) \in \RR^d$ and $\mu(e) \in \RR^d$ 
are defined by $\rho(e) = {\bf e}_i + {\bf e}_j$ and
$\mu(e) = {\bf e}_i - {\bf e}_j$.
Let $A_G \in \ZZ^{d \times n}$ denote the matrix 
with column vectors $\rho(e_1), \ldots, \rho(e_n)$ 
and $A_{\diG} \in \ZZ^{d \times n}$
the matrix with column vectors
$\mu(e_1), \ldots, \mu(e_n)$.
Theorem \ref{finitegraphnormal}
gives a combinatorial characterization on $G$
for which 
$K[A_G^\pm]$ is normal as well as for which  
$I_{A_G^\pm}$ has a squarefree initial ideal.
Theorem \ref{connected bipartite graph}
supplies a quadratic Gr\"obner bases of
the toric ideals of the centrally symmetric configuration
of a bipartite graph any of whose cycles of length $\geq 6$  
has a chord.
Finally, we conclude this paper with several examples of
centrally symmetric configurations of the incidence matrices
of nonbipartite graphs.

\section{The index $[\ZZ^d : \ZZ A]$}

Let $A = [{\bf a}_1 \ \cdots \ {\bf a}_n] \in \ZZ^{d \times n}$ be a 
matrix of rank $d$
and let $\ZZ A$ denote the abelian subgroup of $\ZZ^d$
generated by the column vectors of $A$, i.e.,
$$
\ZZ A := 
\left\{
\left.
\sum_{i=1}^n z_i {\bf a}_i \ \right| \ z_i \in \ZZ, i= 1,2,\ldots,n
\right\}.
$$
Let $[B \  | \  O]$ be the {\it Hermite normal form} 
(\cite[p.45]{Sch})
of $A$
which is obtained by a series of 
elementary unimodular column operations:
\begin{itemize}
\item[(i)]
exchanging two columns;
\item[(ii)]
multiplying a column by $-1$;
\item[(iii)]
adding an integral multiple of one column to another column.
\end{itemize}
Here $B$ is a nonsingular, lower triangular,
nonnegative integer matrix, in which
each row has a unique maximum entry, which is located on the main diagonal of $B$.
As stated in \cite[Proof of Corollary 4.1b]{Sch}, we have $\ZZ A = \ZZ B$.
Moreover, since ``the g.c.d.~of maximal minors" is invariant
under elementary unimodular column operations,
the g.c.d.~of maximal minors of $A$ equals to $|B|$.
Since $[\ZZ^d : \ZZ B] = |B|$, we have the following Propositions:

\begin{Proposition}
\label{gcd}
Let $A \in \ZZ^{d \times n}$ be a 
matrix of rank $d$.
Then the index $[\ZZ^d : ZA]$ equals to the g.c.d.~of maximal minors of $A$.
\end{Proposition}

\begin{Proposition}
\label{henkei}
Let $A \in \ZZ^{d \times n}$ be a 
matrix of rank $d$.
Then there exists a nonsingular matrix $B$ such that
$A' = B^{-1} A $ is an integer matrix satisfying
$\ZZ A' = \ZZ^d$.
\end{Proposition}

Since the centrally symmetric configuration $A^\pm$
of a matrix $A \in \ZZ^{d \times n}$ is
brought into
\[
\left[
\begin{array}{c|ccc|ccc}
0      &   &   &    & &   & \\
\vdots &   & A &    & & O & \\
0      &   &   &    & &   & \\
\hline
1      & 0 & \cdots & 0 & 0 & \cdots & 0
\end{array}
\right]
\]
by a series of 
elementary unimodular column operations,
we have the following:

\begin{Proposition}
\label{indexcsc}
Let $A \in \ZZ^{d \times n}$ be a 
matrix of rank $d$.
Then the index $[\ZZ^{d+1} : \ZZ A^\pm]$ equals to 
the index $[\ZZ^d : \ZZ A]$.
In particular, $\ZZ A^\pm = \ZZ^{d+1}$ if and only if $\ZZ A = \ZZ^d$.
\end{Proposition}

An integer matrix $A \in \ZZ^{d \times n}$ 
of rank $d$ is called {\it unimodular} if 
all nonzero maximal minors of $A$ have the same absolute value.
Let $\delta(A)$ 
denote the absolute value of a nonzero maximal minor
of a unimodular matrix $A$.
For unimodular matrices, Propositions \ref{gcd} and \ref{henkei} 
are

\begin{Corollary}
\label{gcdunimodular}
Let $A \in \ZZ^{d \times n}$ be a unimodular
matrix of rank $d$.
Then the index $[\ZZ^d : \ZZ A]$ equals to $\delta(A)$.
In particular, $\ZZ A = \ZZ^d$ if and only if $\delta(A)=1$.
\end{Corollary}

\begin{Corollary}
\label{henkeiunimodular}
Let $A \in \ZZ^{d \times n}$ be a unimodular
matrix of rank $d$.
Then there exists a nonsingular matrix $B$ such that
$A' = B^{-1} A $ is a unimodular matrix of $\delta(A') =1$. 
\end{Corollary}

\section{Centrally symmetric configurations of unimodular matrices}
Two fundamental results (Theorems \ref{unimodularinitial} and \ref{unimodularGorenstein})
on centrally symmetric configurations 
of unimodular matrices will be established. 

First of all, remark that 
the centrally symmetric configuration $A^\pm$ of a matrix $A$
is NOT unimodular even if $A$ is unimodular.

\begin{Proposition}
\label{cscisnotuni}
Let $A \in \ZZ^{d \times n}$ be a matrix of rank $d$.
Then $A^\pm \in \ZZ^{(d+1) \times (2n+1)}$ is not unimodular.
\end{Proposition}

\begin{proof}
Let $A =
\begin{bmatrix}
{\bf a}_1 & \cdots &{\bf a}_n
\end{bmatrix}
 \in \ZZ^{d \times n}$.
Since ${\rm rank} (A)=d$, there exists a 
nonsingular submatrix
$
A'=
\begin{bmatrix}
{\bf a}_{i_1} &  \cdots & {\bf a}_{i_d} 
\end{bmatrix}
$
of $A$.
Then we have
\begin{eqnarray*}
\begin{vmatrix}
{\bf 0} & {\bf a}_{i_1} & \cdots & {\bf a}_{i_d} \\
1 & 1 & \cdots & 1
\end{vmatrix}
&=&
(-1)^d \ |A'|,\\
\begin{vmatrix}
-{\bf a}_{i_1} & {\bf a}_{i_1}  & \cdots & {\bf a}_{i_d} \\
1 & 1 & \cdots & 1
\end{vmatrix}
&=&
(-1)^d \  2 \ |A'|.
\end{eqnarray*}
Since both of them are nonzero maximal minors of $A^\pm$, 
$A^\pm$ is not unimodular.
\end{proof}

\begin{Example}
{\rm
The centrally symmetric configuration
$$
A^\pm = 
\begin{bmatrix}
0 & 1 & 0 & -1 & 0\\
0 & 0 & 1 & 0 & -1\\
1& 1& 1& 1& 1
\end{bmatrix}
$$ 
of the (totally) unimodular matrix
$A = 
\begin{bmatrix}
1 & 0\\
 0 & 1
\end{bmatrix}
\in \ZZ^{2 \times 2}$ is
not unimodular.
}
\end{Example}

It follows easily that

\begin{Lemma}
\label{interior}
Let $A \in \ZZ^{d \times n}$ be arbitrary.
Then the dimension of
${\rm Conv}(A^\pm) \subset \RR^{d + 1}$
is ${\rm rank}(A)$ and
$[0,\ldots,0,1]^\top \in \RR^{d + 1}$ 
belongs to the interior 
of ${\rm Conv}(A^\pm)$.
\end{Lemma}

Let, in general,
$A \in \ZZ^{d \times n}$ be a configuration
of $\RR^d$ and $B \in \ZZ^{d \times m}$ with $m \leq n$ 
a submatrix (or subconfiguration) of $A$.  
We say that $K[B]$ is a {\em combinatorial pure subring} (see \cite{cpure})
of $K[A]$ if there exists a face $F$ of ${\rm Conv}(A)$ such that
$B = F \cap A$.

\begin{Lemma}
\label{cpuresub}
Suppose that $A \in \ZZ^{d \times n}$ is a configuration.
Then $K[A]$ is a combinatorial pure subring of $K[A^\pm]$.
\end{Lemma}

Recall that $K[A]$ is {\it normal} if and only if 
${\mathbb Z}_{\geq 0} A =
{\mathbb Z} A \cap {\mathbb Q}_{\geq 0} A$ (\cite[Proposition 13.5]{Stu}).
Here ${\mathbb Z}_{\geq 0}$ (resp. ${\mathbb Q}_{\geq 0}$)
is the set of nonnegative integers (resp. nonnegative
rational numbers).
It is known \cite[Proposition 1.2]{cpure} that
if $K[B]$ is a combinatorial pure subring of $K[A]$ 
and if $K[A]$ is normal, then $K[B]$ is normal.

\begin{Corollary}
\label{cpuresubnormalsym}
Suppose that $A \in \ZZ^{d \times n}$ is a configuration
and $K[A^\pm]$ is normal.  Then $K[A]$ is normal.
\end{Corollary}

Example \ref{normalconjecture} stated below shows that the converse of 
Corollary \ref{cpuresubnormalsym} is false.

\begin{Example}
\label{normalconjecture}
{\em
The toric ring $K[A]$ of the configuration
$$A = 
\begin{bmatrix}
2 & 1 & 0\\
0 & 1 & 2
\end{bmatrix}
\in \ZZ^{2 \times 3}$$
is normal.
However, the toric ring $K[A^\pm]$ of the centrally symmetric configuration 
$$
A^{\pm} = 
\begin{bmatrix}
0 & 2 & 1 & 0 & -2 & -1 & 0 \\
0 & 0 & 1 & 2 & 0 & -1 & -2\\
1 & 1 & 1 & 1 & 1 & 1  & 1
\end{bmatrix}
$$
of $A$ is nonnormal.  In fact,
$[1,-1, 1]^\top
\in
{\mathbb Z} A^\pm \cap {\mathbb Q}_{\geq 0} A^\pm $
does not belong to
$
{\mathbb Z}_{\geq 0} A^\pm
$.
}
\end{Example}

The first fundamental result on centrally symmetric configurations 
of unimodular matrices is as follows:

\begin{Theorem}
\label{unimodularinitial}
Let $A \in \ZZ^{d \times n}$ be a unimodular matrix.
Then there exists a reverse lexicographic order $<$ 
on $K[{\bf x}]$ such that the initial ideal
${\rm in}_< (I_{A^\pm})$ of $I_{A^\pm}$ with respect to $<$
is a radical ideal.
\end{Theorem}

\begin{proof}
Let $<$ be a reverse lexicographic order 
on $K[{\bf x}]$ such that
the smallest variable corresponds to 
the column vector $[0,\ldots,0,1]^\top$ of $A^\pm$.
Let $\Delta$ be a pulling triangulation 
\cite[p. 67]{Stu} of ${\rm Conv}(A^\pm)$ 
arising from $<$.
By \cite[Proposition 8.6]{Stu},
the vector
$[0,\ldots,0,1]^\top$
is a vertex of an arbitrary facet (maximal simplex) $\sigma$ 
of $\Delta$.

Let $\ZZ \sigma$ be the abelian subgroup of $\ZZ A^\pm$
spanned by the vertices of $\sigma$.
Since $\sigma$ is a facet of $\Delta$, 
the index $[\ZZ A^\pm : \ZZ \sigma]$ is finite.
The index $[\ZZ A^\pm : \ZZ \sigma]$
is called the {\em normalized volume} of $\sigma$. 
We say that $\Delta$ is {\em unimodular} if 
the normalized volume of each facet $\sigma$ of $\Delta$ 
is equal to $1$. 
Recall that $\Delta$ is unimodular if and only if 
${\rm in}_< (I_{A^\pm})$ is a radical ideal
\cite[Corollary 8.9]{Stu}.
 
Now, our work is to show 
that $\Delta$ is a unimodular triangulation.
Let ${\bf a}_1,\ldots,{\bf a}_n$ be the column vectors 
of $A$ and let
$$
\begin{bmatrix}
{\bf 0}\\
1
\end{bmatrix},
\begin{bmatrix}
\varepsilon_1 {\bf a}_{i_1}\\
1
\end{bmatrix},
\ldots,
\begin{bmatrix}
\varepsilon_d {\bf a}_{i_d}\\
1
\end{bmatrix}
\ \ \ \ (\varepsilon_i \in \{1,-1\})
$$ 
the vertices of $\sigma$.
One has
$$
\begin{vmatrix}
{\bf 0} & \varepsilon_1 {\bf a}_{i_1} & \cdots & \varepsilon_d {\bf a}_{i_d}\\
1 & 1 & \cdots & 1
\end{vmatrix}
=
(-1)^d
\begin{vmatrix}
\varepsilon_1 {\bf a}_{i_1} & \cdots & \varepsilon_d {\bf a}_{i_d}
\end{vmatrix}\\
=
(-1)^d \varepsilon_1 \cdots \varepsilon_{d}
\begin{vmatrix}
 {\bf a}_{i_1} & \cdots & {\bf a}_{i_d}
\end{vmatrix}.
$$
Since $\sigma$ is a simplex, the above determinant cannot be zero.
Hence its absolute value is equal to $\delta(A)$.
Thanks to Proposition \ref{indexcsc} and Corollary \ref{gcdunimodular}, we have
$$
[\ZZ^{d+1} : \ZZ \sigma] = \delta(A) =  [\ZZ^d : \ZZ A]=  [\ZZ^{d+1} : \ZZ A^\pm].
$$
Consequently, one has
$
\ZZ A^\pm 
=
\ZZ
\sigma
$.
In other words,
the normalized volume of $\sigma$ is equal to $1$.
Hence $\Delta$ is a unimodular triangulation 
of ${\rm Conv}(A^\pm)$, as required.
\end{proof}

By virtue of \cite[Proposition 13.15]{Stu}, it follows that

\begin{Corollary}
\label{unimodularnormal}
If $A \in \ZZ^{d \times n}$ is unimodular, then
the toric ring $K[A^\pm]$ of $A^\pm$ is normal.
\end{Corollary}

\begin{Example}
{\rm
Let $A \in \ZZ^{4 \times 5}$ be the configuration
$$
A = 
\begin{bmatrix}
0 & 1 & 1 & 0 & -1\\
0 & 1 & 0 & 1 & -1\\
0 & 0 & 1 & 1 & -1\\
1 & 1 & 1 & 1 & 1 
\end{bmatrix}
$$
and let ${\mathcal P} = {\rm Conv}(A)$.
Then ${\mathcal P}$ is a tetrahedron with ${\mathcal P} \cap \ZZ^4 = A$
and
$[0,0,0,1]^\top$
is a unique integer point 
belonging to the interior of 
 ${\mathcal P}$.
Each facet ${\mathcal F}$ of ${\mathcal P}$ has a trivial unimodular triangulation
since ${\mathcal F}$ is a 2-simplex.
However, $K[A]$ is not normal since $[1,1, 1,2]^\top
\in
{\mathbb Z} A \cap {\mathbb Q}_{\geq 0} A$
does not belong to
$
{\mathbb Z}_{\geq 0} A
$.
Hence, in particular, ${\mathcal P}$ has no unimodular triangulations.
}
\end{Example}

We now turn to the problem of finding Gorenstein toric rings
which arise from centrally symmetric configurations of 
integer matrices.

Let ${\mathcal P} \subset \RR^d$ be an arbitrary convex polytope
of dimension $d$ such that the origin of $\RR^d$ belongs to 
the interior of ${\mathcal P}$.  Then the {\em dual polytope} of 
${\mathcal P} \subset \RR^d$ is defined to be the
convex polytope ${\mathcal P}^\star \subset \RR^d$ 
which consists of those $x \in \RR^d$ with 
$\langle x, y \rangle \leq 1$ 
for all $y \in {\mathcal P}$.  Here $\langle x, y \rangle$
is the usual inner product of $\RR^d$.
One has $({\mathcal P}^\star)^\star = {\mathcal P}$.

An integral convex polytope ${\mathcal P} \subset \RR^d$
of dimension $d$ is called a {\em Fano polytope} 
if the origin of $\RR^d$ is a unique integer point 
belonging to the interior of ${\mathcal P}$.
We say that a Fano polytope ${\mathcal P} \subset \RR^d$
is {\em Gorenstein} if the vertices of 
the dual polytope ${\mathcal P}^\star$
of ${\mathcal P}$ have integer coordinates.

\begin{Remark}
\label{remarkdual}
{\em 
If ${\mathcal P} \subset \RR^d$ is an arbitrary convex polytope
of dimension $d$ such that the origin of $\RR^d$ belongs to 
the interior of ${\mathcal P}$
and if the dual polytope ${\mathcal P}^\star$
is integral, then the origin of $\RR^d$ is a unique integer 
point belonging to the interior of ${\mathcal P}$.
}  
\end{Remark}

Let ${\mathcal P} \subset \RR^N$ be an integral convex polytope
of dimension $d$ and ${\mathcal A} \subset \RR^N$ the affine 
subspace spanned by ${\mathcal P}$.
One has an invertible affine transformation 
$\psi : {\mathcal A} \rightarrow \RR^d$
with $\psi({\mathcal A} \cap \ZZ^N) = \ZZ^d$. 
It follows that $\psi({\mathcal P}) \subset \RR^d$
is an integral convex polytope of dimension $d$.
We say that $\psi({\mathcal P})$ is a {\em standard form}
of ${\mathcal P}$.
In general, by abuse of terminology,
we say that an integral convex polytope 
is a Gorenstein Fano polytope if one of its standard forms 
is a Gorenstein Fano polytope.

\begin{Lemma}
\label{uniGorlemma}
Suppose that $A \in \ZZ^{d \times n}$ is a unimodular matrix
of rank $d$ with $\delta(A) = 1$.
Then the integral convex polytope
${\rm Conv}(A^{\pm}) \subset \RR^{d + 1}$
of dimension $d$ is a Gorenstein Fano polytope.
\end{Lemma}

\begin{proof}
Let ${\mathcal A}$ be the affine subspace of $\RR^{d+1}$
spanned by ${\rm Conv}(A^{\pm})$.  In other words,
${\mathcal A}$ is the hyperplane of $\RR^{d+1}$
defined by the equation $z_{d+1} = 1$.
Let $\psi : {\mathcal A} \rightarrow \RR^d$
be the invertible affine transformation
defined by 
\[
\psi(z_1, \ldots, z_d, z_{d+1})
= (z_1, \ldots, z_d).
\]
It follows from Lemma \ref{interior} that
the standard form $\psi({\rm Conv}(A^{\pm})) \subset \RR^d$ 
contains the origin of $\RR^d$ in its interior.
Let ${\mathcal F}$ be a facet of $\psi({\rm Conv}(A^{\pm}))$
and $\sigma \subset {\mathcal F}$ an arbitrary $(d - 1)$-simplex.   
Now, by virtue of the proof of Theorem \ref{unimodularinitial},
the vertices of $\sigma$ is a $\ZZ$-basis of $\ZZ^{d}$.  
In particular
the equation of the facet ${\mathcal F}$
is of the form $a_1z_1 + \cdots + a_{d}z_{d} = 1$ 
with each $a_i \in \ZZ$.  In other words,
the dual polytope of $\psi({\rm Conv}(A^{\pm}))$ is integral.
Hence $\psi({\rm Conv}(A^{\pm}))$ is a Gorenstein Fano polytope.
\end{proof}

\begin{Example}
\label{GorDeltaAis2}
{\em
Let $A \in \ZZ^{2 \times 2}$ be the unimodular matrix
$$
A = 
\begin{bmatrix}
1 & 1\\
1 &-1  
\end{bmatrix}
$$
with $\delta(A)=2$.
Then the convex hull of  $A^\pm$
is a Gorenstein Fano polytope.
}
\end{Example}

\begin{Example}
\label{GorEX}
{\em
Let $A \in \ZZ^{3 \times3}$ and $B \in \ZZ^{3 \times 3}$ 
be the unimodular matrices  
\[
A = 
\left[
\begin{array}{ccc}
1 & 1 & 0 \\
1 & 0 & 1 \\
0 & 1 & 1  
\end{array}
\right], \, \, \, \, \, 
B = 
\left[
\begin{array}{ccc}
1 & 1 & 0 \\
1 & 0 & 1 \\
1 & 1 & 1 
\end{array}
\right]
\]
with $\delta(A) = 2$ and $\delta(B) = 1$.
One has  
\[
I_{A^\pm} = I_{B^\pm} = 
\langle
x_1^2 - x_2x_5, x_1^2 - x_3x_6, x_1^2 -  x_4x_7
\rangle.
\]
Lemma \ref{uniGorlemma} says  
that ${\rm Conv}(B^\pm)$ is a Gorenstein Fano polytope.
However, ${\rm Conv}(A^\pm)$
cannot be a Gorenstein Fano polytope. 
}
\end{Example}

Let ${\mathcal P} \subset \RR^d$ be an integral convex polytope.
For $N = 1,2,\ldots$, we define
$$
N {\mathcal P} = 
\{
N \alpha \in \RR^d
\ | \ 
\alpha \in  {\mathcal P}
\}.
$$

\begin{Lemma}
\label{basiclemma}
Let $A \in \ZZ^{d \times n}$ be a configuration with $\ZZ A = \ZZ^d$
and suppose that $K[A]$ is normal.  
Then, for each $\alpha \in N {\rm Conv} (A) \cap \ZZ^d$,
there exist $\beta_1, \ldots, \beta_N \in  {\rm Conv} (A) \cap \ZZ^d$ such that
$ \alpha = \beta_1 + \cdots + \beta_N$.
\end{Lemma}

\begin{proof}
Let $A=[ {\bf a}_1 \ \cdots \ {\bf a}_n]$ and suppose that
$\alpha$ belongs to $N {\rm Conv} (A) \cap \ZZ^d$.
Since $\alpha$ belongs to $N {\rm Conv} (A)$, 
one has $\alpha = \sum_{i=1}^n r_i {\bf a}_i $, 
where $0 \leq r_i \in \QQ$ and where
$\sum_{i=1}^n r_i = N$.
In particular $\alpha $ belongs to 
$\ZZ^d \cap \QQ_{\geq 0} A$.
Since 
$\ZZ A = \ZZ^d$ and since $K[A]$ is normal,
one has
$
\ZZ^d \cap \QQ_{\geq 0} A 
=
\ZZ A \cap \QQ_{\geq 0} A 
=
\ZZ_{\geq 0} A.
$
Thus 
$
\alpha = \sum_{i=1}^n r_i {\bf a}_i   = \sum_{i=1}^n z_i {\bf a}_i 
$
with each
$z_i \in \ZZ_{\geq 0}$.
Since $A$ is a configuration, it follows easily that  
$\sum_{i=1}^n z_i = N$, as required.
\end{proof}

We now come to the second fundamental result on centrally symmetric 
configurations of unimodular matrices.
 
\begin{Theorem}
\label{unimodularGorenstein}
Suppose that $A \in \ZZ^{d \times n}$ is a unimodular matrix
of rank $d$.
Then the toric ring $K[A^\pm]$ of the centrally symmetric configuration 
$A^{\pm}$ is Gorenstein.
\end{Theorem}

\begin{proof}
By virtue of Corollary \ref{henkeiunimodular},
we may assume that $\delta(A) = 1$.
Lemma \ref{uniGorlemma} says that the integral convex  
${\rm Conv}(A^{\pm}) \subset \RR^{d+1}$ is 
a Gorenstein Fano polytope.  Hence 
\cite[Corollary (1.2)]{DeNegriHibi} guarantees that
the Ehrhart ring \cite[p. 97]{Hibi} 
of ${\rm Conv}(A^{\pm})$ is Gorenstein.
Since $\delta(A) = 1$, one has 
$\ZZ A^\pm = \ZZ^{d+1}$.
Moreover, the toric ring $K[A^\pm]$ is normal.
Thus, by using Lemma \ref{basiclemma}, it follows that $K[A^\pm]$ coincides with
the Ehrhart ring of ${\rm Conv}(A^{\pm})$.
Hence $K[A^\pm]$ is Gorenstein, as desired.   
\end{proof}

\begin{Example}
\label{GorensteinConverse}
{\em
Let $A \in \ZZ^{2 \times 4}$ be the matrix
$$
A = 
\begin{bmatrix}
0 & 1 & 1 & 1\\
1 & 0 & 1 &-1
\end{bmatrix}
$$
which is not unimodular.
The toric ring $K[A^\pm]$ of 
$$
A^\pm = 
\begin{bmatrix}
0 & 0 & 1 & 1 & 1 &  0 &-1 &-1 &-1\\ 
0 & 1 & 0 & 1 &-1 &-1 &  0 &-1 & 1\\
1 & 1 &  1 &  1 &  1 & 1 &  1 &  1 &  1 
\end{bmatrix}
$$
is normal and Gorenstein.
}
\end{Example}

\section{Centrally symmetric configurations of finite graphs}

Let $G$ be a finite connected graph on the vertex set
$[d] = \{ 1, \ldots, d \}$ and suppose that $G$ possesses
no loop and no multiple edge.
Let $E(G) = \{ e_1, \ldots, e_n \}$ denote
the edge set of $G$.  Let ${\bf e}_1, \ldots, {\bf e}_d$
stand for the canonical unit coordinate vector of $\RR^d$.
If $e = \{ i, j \}$ is an edge of $G$ with $i < j$, then
the column vectors $\rho(e) \in \RR^d$ and $\mu(e) \in \RR^d$ 
are defined by $\rho(e) = {\bf e}_i + {\bf e}_j$ and
$\mu(e) = {\bf e}_i - {\bf e}_j$.
Let $A_G \in \ZZ^{d \times n}$ denote the matrix 
with column vectors $\rho(e_1), \ldots, \rho(e_n)$ 
and $A_{\diG} \in \ZZ^{d \times n}$
the matrix with column vectors
$\mu(e_1), \ldots, \mu(e_n)$.
Thus $A_G$ is a configuration of $\RR^d$.
However, $A_{\diG}$ is not necessarily a configuration of $\RR^d$.

The $(0, 1)$-polytope
${\rm Conv}(A_{G}) \subset \RR^{d}$ is called 
the {\em edge polytope} \cite{edgepolytope} of $G$.
We say that 
${\rm Conv}({A_{\diG}^\pm}) \subset \RR^{d+1}$ 
is the {\em symmetric edge polytope} of $G$.
The symmetric edge polytope
${\rm Conv}({A_{\diG}^\pm}) \subset \RR^{d+1}$ is 
of dimension $d - 1$.
The dimension of ${\rm Conv}(A_G)$ is given in \cite[Proposition 1.3]{edgepolytope}.
Note that $\dim({\rm Conv}({A_{G}^\pm})) = 1+ \dim( {\rm Conv}(A_G) )$. 
If $G$ is nonbipartite graph,
then the dimension of 
${\rm Conv}({A_{G}^\pm}) \subset \RR^{d+1}$ is $d$.
If $G$ is bipartite, 
then the dimension of 
${\rm Conv}({A_{G}^\pm}) \subset \RR^{d+1}$ is $d - 1$.

Recall that an integer matrix is {\em totally unimodular}
(\cite[p.266]{Sch}) if every square submatrix 
has determinant $0$, $1$, or $- 1$.  In particular
every entry of a totally unimodular matrix belongs to
$\{ 0, 1, - 1 \}$.  
It is known that $A_{\diG}$ is totally unimodular.
In addition,  $A_G$ is totally unimodular if and only if $G$ is a bipartite graph.

It follows from Theorem \ref{unimodularinitial}
and Corollary \ref{unimodularnormal} that

\begin{Corollary}
\label{diG}
Let $G$ be a finite connected graph.
Then there exists a reverse lexicographic order $<$ 
on $K[{\bf x}]$ such that
${\rm in}_<(I_{A_{\diG}^\pm})$ is a radical ideal.
Thus in particular $K[A_{\diG}^\pm]$ is normal.
\end{Corollary}

On the other hand, it is known \cite[Corollary 2.3]{edgepolytope}
that

\begin{Proposition}
\label{graphnormalold}
{\rm
Let $G$ be a finite connected graph.
Then $K[A_G]$ is normal if and only if,
for each two odd cycles $C_1$ and $C_2$ of $G$ having no common vertex,
there exists an edge of $G$ which joins a vertex of $C_1$ and a vertex of $C_2$.
}
\end{Proposition}

We now discuss when 
the toric ring $K[A_G^\pm]$ is normal.

\begin{Theorem}
\label{finitegraphnormal}
Let $G$ be a finite connected graph.
Then the following conditions are equivalent:
\begin{itemize}
\item[(i)]
$K[A_G^\pm]$ is normal;
\item[(ii)]
$I_{A_G^\pm}$ has a squarefree initial ideal;
\item[(iii)]
$A_G$ is a unimodular matrix
(by deleting a redundant row if $G$ is bipartite);
\item[(iv)]
Any two odd cycles of $G$ possess a common vertex.
\end{itemize}
\end{Theorem}

\begin{proof}
First, (ii) $\Longrightarrow $ (i) 
is known (\cite[Proposition 13.15]{Stu}). 
Second, (iii) $\Longleftrightarrow $ (iv)
is discussed in, e.g., \cite{minorsgraph}. 
Third, 
(iii) $\Longrightarrow $ (ii)
follows from Theorem \ref{unimodularinitial}.

Now, in order to show
(i) $\Longrightarrow $ (iv),
suppose that $K[A_G^\pm]$ is normal
and that there exist two odd cycles 
$C_1 = (i_1,\ldots,i_r)$ and 
$C_2=(j_1,\ldots,j_s)$ 
of $G$ having no common vertex.
Thanks to Corollary \ref{cpuresubnormalsym}
together with Proposition \ref{graphnormalold},
there exists an edge $e$ of $G$ which joins a vertex of $C_1$ and a vertex of $C_2$.
Let, say, $e = \{i_r, j_s\}$.
Then
\begin{eqnarray*}
\alpha &=& 
\sum_{k=1}^{\frac{r-1}{2}}
({\bf e}_{i_{2k-1}} + {\bf e}_{i_{2k}} +{\bf e}_{d+1} )
+
\sum_{\ell=1}^{\frac{s-1}{2}}
(-{\bf e}_{j_{2 \ell}} - {\bf e}_{j_{2 \ell+1}} +{\bf e}_{d+1} )
\\
 & & 
+
(-{\bf e}_{j_s} - {\bf e}_{j_1} +{\bf e}_{d+1} )
+
({\bf e}_{i_r} +{\bf e}_{j_s} +{\bf e}_{d+1} )- {\bf e}_{d+1}\\
&=&
{\bf e}_{i_1} + {\bf e}_{i_2} +\cdots + {\bf e}_{i_r}
-{\bf e}_{j_1} - {\bf e}_{j_2} -\cdots - {\bf e}_{j_s}
+\frac{r+s}{2} {\bf e}_{d+1} 
\end{eqnarray*}
belongs to ${\mathbb Z} A_G^\pm$.
Let $i_{r+1} = i_1$ and $j_{s+1} = j_1$.
It then follows that
$$
\alpha =
\sum_{k=1}^{r}
\frac{1}{2}
({\bf e}_{i_{k}} + {\bf e}_{i_{k+1}} +{\bf e}_{d+1} )
+
\sum_{\ell=1}^{s}
\frac{1}{2}
(-{\bf e}_{j_{\ell}} - {\bf e}_{j_{\ell+1}} +{\bf e}_{d+1} )
\in
{\mathbb Q}_{\geq 0} A_G^\pm.
$$
Since $K[A_G^\pm]$ is normal, 
$\alpha$ belongs to ${\mathbb Z}_{\geq 0} A_G^\pm$.
Thus 
$$
\alpha = 
z_0 
\begin{bmatrix} {\bf 0} \\ 1 \end{bmatrix}
+
\sum_{i=1}^n z_i 
\begin{bmatrix} \rho(e_i) \\ 1 \end{bmatrix}
+
\sum_{j=1}^n z_j'
\begin{bmatrix} - \rho(e_j)  \\ 1 \end{bmatrix}
$$
where $0 \leq z_i , z_j' \in \ZZ$.
Since both $r$ and $s$ are odd and
$\{i_1,\ldots ,i_r\} \cap \{j_1,\ldots,j_s\} = \emptyset$,
it follows that
$
\frac{r+1}{2} \leq \sum_{i=1}^n z_i
$
and
$
\frac{s+1}{2} \leq \sum_{j=1}^n z_j'.
$
Hence
$$
\frac{r+1}{2} + \frac{s+1}{2}
\leq 
z_0+ \sum_{i=1}^n z_i+  \sum_{j=1}^n z_j'
=
\frac{r+s}{2},
$$
a contradiction.
\end{proof}

It follows from the theory of totally unimodular matrices
that the toric ring $K[A_{\diG}^\pm]$ of 
the centrally symmetric configuration  
$A_{\diG}^\pm$ of
a finite connected graph $G$ is Gorenstein.
Moreover, if $G$ is bipartite, then
$K[A_{G}^\pm]$ is Gorenstein.

\begin{Theorem}
\label{nonbipartiteGor}
Let $G$ be a finite connected graph
and suppose that any two odd cycles
of $G$ possesses a common vertex.  Then the toric ring
$K[A_G^\pm]$ of the centrally symmetric configuration $A_G^\pm$ 
is Gorenstein.
\end{Theorem}

\begin{proof}
By virtue of the equivalence of (iii) and (iv) of 
Theorem \ref{finitegraphnormal},
it follows that the matrix $A_G$ is unimodular.
Hence Theorem \ref{unimodularGorenstein}
guarantees that $K[A_G^\pm]$
is Gorenstein, 
as desired.
\end{proof}

\begin{Example}
\label{Hilbertfunction}
{\em 
Let $W_d$ be the wheel graph on $[d]$.
For example, $W_6$ is as follows:

\begin{figure}[h]
\begin{center}
\includegraphics[width=3cm]{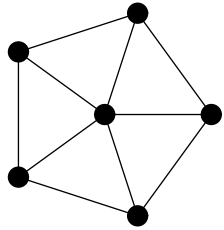}
\caption{The wheel graph $W_6$}
\end{center}
\end{figure}

\noindent
Thanks to Theorems \ref{finitegraphnormal} and \ref{nonbipartiteGor},
$K[ A_{W_d}^\pm ]$ is normal and Gorenstein.
One can compute the Hilbert series
$$
H ( K[ A_{W_d}^\pm ], \lambda) =
\sum_{j=0}^\infty \dim_K \left( K[ A_{W_d}^\pm ]  \right)_j \lambda^j
$$
of $K[ A_{W_d}^\pm ]$
\begin{table}[h]
\begin{center}
\begin{tabular}{|c|c|}
\hline
$d$ &  $ (1-\lambda)^{d+1} H ( K[ A_{W_d}^\pm ], \lambda)$ \\
\hline
4 & $1 + 8\lambda + 14\lambda^2 + 8\lambda^3 + \lambda^4 $ \\
\hline
5 & $ 1 + 11\lambda + 32\lambda^2 + 32\lambda^3 + 11\lambda^4 + \lambda^5 $ \\
\hline
6 & $1 + 14\lambda + 65\lambda^2 + 104\lambda^3 + 65\lambda^4 + 14\lambda^5 + \lambda^6$ \\
\hline
7 & $1 + 17\lambda + 105\lambda^2 + 249\lambda^3 + 249\lambda^4 + 105\lambda^5 + 17\lambda^6 + \lambda^7$\\
\hline
\end{tabular}

\bigskip

\caption{The Hilbert series of the toric ring of the wheel graph $W_d$}
\end{center}
\end{table}
by the software {\tt CoCoA} \cite{cocoa}.
}
\end{Example}

\begin{Example}
{\em 
Let $G$ be the graph on the vertex set $\{1,\ldots,6\}$ with 
the edge set
$$
E(G) =
\{
\{1 ,2 \},
\{2 ,3 \},
\{1 ,3 \},
\{3 ,4 \},
\{ 4, 5\},
\{5 ,6 \},
\{ 4, 6\}
\}.
$$
By virtue of Theorem \ref{finitegraphnormal},
$K[ A_G^\pm ]$ is not normal.
On the other hand,  by {\tt CoCoA}, one can check that
$K[A_G^\pm]$ is Cohen--Macaulay and not Gorenstein.
}
\end{Example}

\begin{Conjecture}
{\em
Let $G$ be a fnite connected graph.
Then, 
$K[A_G^\pm]$ is Gorenstein if and only if 
any two odd cycles
of $G$ possesses a common vertex.
}
\end{Conjecture}

\begin{Question}
{\em
Let $G$ be a fnite connected graph such that $A_G$ is not unimodular.
What is the necessary and sufficient condition
for $K[A_G^\pm]$ to be Cohen--Macaulay?
}
\end{Question}

\section{Centrally symmetric configurations of bipartite graphs}

In this section, we study toric ideals of 
centrally symmetric configurations arising from bipartite graphs.

\begin{Proposition}
Let $G$ be a bipartite graph.
Then,
we have $I_{A_G^\pm} = I_{A_{\diG}^\pm}$.
\end{Proposition}

\begin{proof}
If $G$ is bipartite, 
the rows of $A_G$ equals to
the rows of $A_{\diG}$ up to $- 1$ multiples.
Hence 
the rows of $A_G^\pm$ equals to
the rows of $A_{\diG}^\pm$ up to $- 1$ multiples.
Thus we have ${\rm Ker} \ A_G^\pm = {\rm Ker} \ A_{\diG}^\pm$.
By \cite[Corollary 4.3]{Stu}, it follows that $I_{A_G^\pm} = I_{A_{\diG}^\pm}$.
\end{proof}

By Lemma \ref{cpuresub}, since $A_G$ is a configuration,
$K[A_G]$ is a combinatorial pure subring of $K[A_G^\pm]$.
It is known \cite{koszulbipartite} that

\begin{Proposition}
\label{koszul}
Let $G$ be a connected bipartite graph.
Then the following conditions are equivalent:
\begin{itemize}
\item[(i)]
Every cycle of length $\geq 6$ in $G$ has a chord;
\item[(ii)]
$I_{A_G}$ has a quadratic Gr\"obner basis;
\item[(iii)]
$K[A_G]$ is Koszul;
\item[(iv)]
$I_{A_G}$ is generated by quadratic binomials.
\end{itemize}
\end{Proposition}

Thus we have the following.

\begin{Corollary}
Let $G$ be a connected bipartite graph.
If $G$ has an cycle of length $\geq 6$ 
having no chord,
then $I_{A_G^\pm}$ 
{\rm (}$=I_{A_{\diG}^\pm}${\rm )} is not generated by quadratic binomials.
\end{Corollary}

Let $G$ be a connected bipartite graph
on the vertex set $\{1,\ldots,p\} \cup \{1',\ldots,q'\}$.
Suppose that every cycle of $G$ of length $\geq 6$ 
has a chord.
Then, by the same argument as in \cite{koszulbipartite},
we may assume that 
$G$ satisfies the condition $(*)$:
$$
\{i,\ell'\} ,\{j,k'\} ,\{j,\ell'\} \in E(G) 
\ \ 
\Longrightarrow 
\ \ 
\{i,k'\} \in  E(G)
\eqno{(*)}
$$
for each $1 \leq i < j \leq p$, $1 \leq k < \ell \leq q$.
Let ${\mathcal R} = K[s_1^{\pm 1},\ldots,s_p^{\pm 1},t_1^{\pm 1},\ldots,t_q^{\pm 1},u]$
and
$
K[A_{\diG}^\pm]
=
K[\{u\} \cup \{s_i t_j^{-1} u \ | \ \{i,j' \} \in E(G)\} \cup \{s_i^{-1} t_j u \ | \ \{i,j' \} \in E(G)\}]
\subset {\mathcal R}.$
Let ${\mathcal S} = K[
\{ x_{i j} \}_{ \{i,j'\} \in E(G)} 
\cup
\{ y_{i j} \}_{ \{i,j'\} \in E(G)} 
\cup
\{z\}
]$ be a polynomial ring over $K$.
Then 
the toric ideal $I_{A_{\diG}^\pm}$
($=I_{A_G^\pm}$)
is the kernel of surjective homomorphism
$\varphi : {\mathcal S} \rightarrow K[A_{\diG}^\pm]$
defined by $\varphi (x_{ij}) = s_i t_j^{-1} u$, $\varphi (y_{ij}) = s_i^{-1} t_j u$
 and $\varphi (z) = u $.
Let $<$ denote the reverse lexicographic order on ${\mathcal S}$
induced by the ordering
$$
z < y_{11} < x_{11} < y_{12} < x_{12} < \cdots < y_{1q} < x_{1q} < 
y_{21} < x_{21} < \cdots < y_{pq} < x_{pq}.  
$$

\begin{Theorem}
\label{connected bipartite graph}
Let $G$ be a connected bipartite graph.
Suppose that
every cycle in $G$ of length $\geq 6$ 
has a chord.
Let ${\mathcal G}$ be the set consists of the following binomials:
\begin{eqnarray*}
x_{ik} y_{ik} -z^2 &  & \{i,k'\} \in E(G)\\
x_{i \ell} x_{j k} - x_{i k} x_{j \ell} &  & \{i,k'\},\{i,\ell'\},\{j,k'\},\{j,\ell'\} \in E(G),\ \ 
i < j  , \ \   k < \ell \\
x_{i \ell} y_{j \ell} - x_{i k} y_{j k} &  & \{i,k'\},\{i,\ell'\},\{j,k'\},\{j,\ell'\} \in E(G),\ \ 
i < j  , \ \   k < \ell \\
y_{j \ell} x_{j k} - x_{i k} y_{i \ell} &  & \{i,k'\},\{i,\ell'\},\{j,k'\},\{j,\ell'\} \in E(G),\ \ 
i < j  , \ \   k < \ell \\
x_{j \ell} y_{j k} - y_{i k} x_{i \ell} &  & \{i,k'\},\{i,\ell'\},\{j,k'\},\{j,\ell'\} \in E(G),\ \ 
i < j  , \ \   k < \ell \\
y_{i \ell} x_{j \ell} - y_{i k} x_{j k} &  & \{i,k'\},\{i,\ell'\},\{j,k'\},\{j,\ell'\} \in E(G),\ \ 
i < j  , \ \   k < \ell \\
y_{i \ell} y_{j k} - y_{i k} y_{j \ell} &  & \{i,k'\},\{i,\ell'\},\{j,k'\},\{j,\ell'\} \in E(G),\ \ 
i < j  , \ \   k < \ell 
\end{eqnarray*}
where the initial monomial of each binomial is the first monomial and squarefree.
Then ${\mathcal G}$ is 
the reduced Gr\"obner basis of $I_{A_{\diG}^\pm}$ with respect to $<$.
\end{Theorem}

\begin{proof}
It is easy to see that ${\mathcal G} \subset I_{A_{\diG}^\pm}$ and that
 the initial monomial of each binomial in ${\mathcal G}$ is the first monomial and squarefree.

Suppose that ${\mathcal G}$ is not the reduced Gr\"obner basis of $I_{A_{\diG}^\pm}$ with respect to $<$.
Then there exists an irreducible binomial $f=m_1-m_2 \in I_{A_{\diG}^\pm} $ such that, for $i = 1,2$, 
$m_i$ is not divisible by
the initial monomial of any binomial in ${\mathcal G}$.

Suppose that the biggest variable appearing in $m_1 m_2$ is $x_{j \ell}$.
We may assume that $m_1$ is divided by $x_{j \ell}$ and 
$m_2 $ is not divided by $x_{j \ell}$.
Since $m_1$ is not divided by $x_{j \ell}y_{j \ell}$, 
$m_1$ is not divided by $y_{j \ell}$.
Let $\varphi (m_1 ) = s_1^{\alpha_1} \cdots s_p^{\alpha_p} t_1^{\beta_1} \cdots t_q^{\beta_q} u^\gamma $.

\bigskip

\noindent
{\bf Case 1.}
$m_1$ is divided by $y_{j k}$ for some $k < \ell$.

Suppose that $\{i,\ell'\} \in E(G)$ for some $i < j$.
Thanks to the condition $(*)$, we have $\{i,k'\} \in E(G)$.
Hence, $m_1$ is divided by the initial monomial of the binomial
$
x_{j \ell} y_{j k} - y_{i k} x_{i \ell}
$
and this is a contradiction.
Thus, $\{i,\ell'\} \notin E(G)$ for all $i < j$.
Then we have $\beta_\ell <0$.
Since $\varphi (m_1) = \varphi (m_2)$,
$m_2$ is divided by $x_{\lambda \ell}$ for some $\lambda < j$.
Then, $\{\lambda,\ell'\} \in E(G)$ and this is a contradiction.

\bigskip

\noindent
{\bf Case 2.}
$m_1$ is divided by $y_{i \ell}$ for some $i < j$.

Similar to Case 1.
Suppose that $\{j,k'\} \in E(G)$ for some $k < \ell$.
Thanks to the condition $(*)$, we have $\{i,k'\} \in E(G)$.
Hence, $m_1$ is divided by the initial monomial of the binomial
$
y_{i \ell} x_{j \ell} - y_{i k} x_{j k}
$
and this is a contradiction.
Thus, $\{j,k'\} \notin E(G)$ for all $k < \ell$.
Then we have $\alpha_j > 0$.
Since $\varphi (m_1) = \varphi (m_2)$,
$m_2$ is divided by $x_{j \mu}$ for some $\mu < \ell$.
Then, $\{j, \mu'\} \in E(G)$ and this is a contradiction.

\bigskip

\noindent
{\bf Case 3.}
$m_1$ is not divided by $y_{j k}$ for all $k < \ell$ and
not divided by $y_{i \ell}$ for all $i < j$.

It then follows that $\alpha_j > 0$ and $\beta_\ell < 0$.
Since $\varphi (m_1) = \varphi (m_2)$,
$m_2$ is divided 
by $x_{i \ell}$ for some $i < j$
and
by $x_{j k}$ for some $k < \ell$.
Thanks to the condition $(*)$, we have $\{i,k'\} \in E(G)$.
Hence, $m_2$ is divided by the initial monomial of the binomial
$
x_{i \ell} x_{j k} - x_{i k} x_{j \ell}$
and this is a contradiction.

\bigskip

On the other hand, if the biggest variable appearing in $m_1 m_2$ is $y_{j \ell}$,
then, by the similar argument (changing the role of $x$ and $y$) as above,
a contradiction arises.
\end{proof}

\begin{Corollary}
\label{symmetricedgeKoszul}
Let $G$ be a connected bipartite graph and let $A = A_{\diG}^\pm$.
Then the following conditions are equivalent:
\begin{itemize}
\item[(i)]
Every cycle of length $\geq 6$ in $G$ has a chord;
\item[(ii)]
$I_A$ has a quadratic Gr\"obner basis;
\item[(iii)]
$K[A]$ is Koszul;
\item[(iv)]
$I_A$ is generated by quadratic binomials.
\end{itemize}
\end{Corollary}

\section{Examples of nonbipartite graphs}

In the previous section, we showed that,
if $G$ is a bipartite graph, then the following conditions are equivalent:
\begin{itemize}
\item[(i)]
$I_{A_G}$ is generated by quadratic binomials;
\item[(ii)]
$I_{A_G^\pm}$ is generated by quadratic binomials;
\item[(iii)]
$I_{A_{\diG}^\pm}$ is generated by quadratic binomials.
\end{itemize}

In this section, we study toric ideals of 
centrally symmetric configurations arising from nonbipartite graphs.
Since $K[A_G]$ is a combinatorial pure subring of $K[A_G^\pm]$,
(ii) $\Longrightarrow $ (i) holds for a nonbipartite graph $G$.
However, if $G$ is not  bipartite, then
none of the other directions hold.
(Computation below is done by {\tt CoCoA}.)

\begin{Example}
{\em
Let $G$ be a graph on the vertex set $\{1,\ldots,6\}$ together with the edge set
$
E(G) = \{
\{1,2\},
\{2 , 3\},
\{3 ,4 \},
\{4 , 5\},
\{5 , 6\},
\{1 , 6\},
\{1 , 3\},\{1 ,5 \},\{2,6\}
\}
.$
Then, 
$I_{A_G}$
has a Gr\"obner basis consists of 3 quadratic binomials.
However, 
neither
$I_{A_G^\pm}$ nor
$I_{A_{\diG}^\pm}$ is generated by quadratic binomials.
Thus neither ``(i) $\Longrightarrow $ (ii)" nor ``(i) $\Longrightarrow $ (iii)" hold.
}
\end{Example}

\begin{Example}
{\em
Let $G$ be a complete 3-partite graph on the vertex set $\{1,2\} \cup \{3,4\}\cup \{5,6\}$.
Then, $I_{A_G^\pm}$ is generated by quadratic binomials.
However, 
$I_{A_{\diG}^\pm}$ is not generated by quadratic binomials.
Thus ``(ii) $\Longrightarrow $ (iii)" does not hold.
}
\end{Example}

\begin{Example}
\label{tie}
{\em
Let $G$ be a graph on the vertex set $\{1,\ldots,5\}$ together with the edge set
$
E(G) = \{
\{1,5\},
\{3,5\},
\{1,3\},
\{2,5 \},
\{4, 5\},
\{2,4\}
\}
.$
Then, 
$I_{A_{\diG}^\pm}$ is generated by quadratic binomials.
However, 
$I_{A_G}= \left< x_1 x_2 x_6  - x_3 x_4 x_5 \right>$ and hence
$I_{A_G^\pm}$ is not generated by quadratic binomials.
Thus neither ``(iii) $\Longrightarrow $ (i)" nor ``(iii) $\Longrightarrow $ (ii)" hold.
}
\end{Example}

\begin{Remark}
{\em
A combinatorial criterion for the toric ideal $I_{A_G}$ of a finite connected graph $G$
to be generated by quadratic binomials is given in \cite[Theorem 1.2]{quadgene}.
}
\end{Remark}

\begin{Proposition}
\label{essbipartite}
Let $G$ be a nonbipartite graph.
Suppose that there exists a vertex $v$ of $G$ such that
any odd cycle of $G$ contains $v$.
Then there exists a bipartite graph $G'$ such that
$
I_{A_G} = I_{A_{G'}} 
$
and
$
I_{A_G^\pm} = I_{A_{G'}^\pm}.
$
\end{Proposition}

\begin{proof}
Let $V(G) = \{1,2,\ldots,d\}$.
Suppose that any odd cycle of $G$ contains the vertex $d$.
Then the induced subgraph $G''$ of $G$ on the vertex set $\{1,2,\ldots,d-1\}$
is a bipartite graph.
(Note that $G''$ is not necessarily connected.)
Let $\{1,2,\ldots,d-1\}=V_1 \cup V_2$ denote a partition of the vertex set of $G''$.
Let $G'$ be a bipartite graph on the vertex set $\{1,2,\ldots,d,d+1\}$ together 
with the edge set
$$ E(G'') \cup \{ \{i, d\} \in E(G) \ | \ i \in V_1  \} 
\cup \{ \{i, d+1\}  \ | \ i \in V_2,  \{i, d\} \in E(G) \}.$$
It then follows that $A_{G'}$ is brought into $A_G$ by
adding the $(d+1)$th row to the $d$th row and deleting the $(d+1)$th row
which is redundant.
\end{proof}

\begin{Example}
{\em
Let $G$ be the nonbipartite graph in Example \ref{tie}
and let $G'$ be a cycle of length 6.
Then $G$ and $G'$ satisfy the condition in Proposition \ref{essbipartite}.
The corresponding matrices are
$$
A_G =
\begin{bmatrix}
1&0&1&0&0&0\\
0&0&0&1&0&1\\
0&1&1&0&0&0\\
0&0&0&0&1&1\\
\hline
1&1&0&1&1&0
\end{bmatrix}
, \ \ 
A_{G'} = 
\begin{bmatrix}
1&0&1&0&0&0\\
0&0&0&1&0&1\\
0&1&1&0&0&0\\
0&0&0&0&1&1\\
\hline
1&0&0&1&0&0\\
0&1&0&0&1&0
\end{bmatrix}.
$$
}
\end{Example}

\begin{Example}
{\em
The wheel graph $W_d$ on $[d]$ satisfies the condition 
in Proposition \ref{essbipartite}
if and only if $d$ is odd.
On the other hand,
it is known \cite[Example 2.1]{quadgene} that $I_{A_{W_6}}$ is generated by quadratic binomials
and has no quadratic Gr\"obner basis.
Thanks to Proposition \ref{koszul}, there exists no bipartite graph $G'$
such that $I_{A_{W_6}} = I_{A_{G'}}$.
}
\end{Example}

\begin{Remark}
{\em
If $G$ is the complete graph on the vertex set $[d]$,
then $A_{\diG}^\pm$ coincides with 
the configuration $M_{A_d}'$ studied in \cite{hosten}.
}
\end{Remark}

\section*{Acknowledgement}
This research was supported by 
the JST
(Japan Science and Technology Agency)
CREST 
(Core Research for Evolutional Science and Technology)
research project
{\em Harmony of Gr\"obner Bases and the Modern Industrial Society}
in the frame of the Mathematic Program
``Alliance for Breakthrough between Mathematics and Sciences.''

\bigskip
\bigskip

\noindent
Hidefumi Ohsugi\\
Department of Mathematics\\
College of Science\\
Rikkyo University\\
Toshima-ku, Tokyo 171-8501, Japan\\
{\tt ohsugi@rikkyo.ac.jp}

\bigskip

\noindent
Takayuki Hibi\\
Department of Pure and Applied Mathematics\\
Graduate School of Information Science and Technology\\
Osaka University\\
Toyonaka, Osaka 560-0043, Japan\\
{\tt hibi@math.sci.osaka-u.ac.jp}
\end{document}